# Comprehensive Modeling of Electric Vehicles in California Demand Response Markets

Bin Wang, Rongxin Yin, and Doug Black

*Abstract*—Electric vehicle (EV) is a significant type of distributed energy resources (DERs), that provide flexibilities to grid operators to achieve a myriad of objectives. This paper presents a comprehensive modeling framework of EVs under multiple real-world demand response (DR) markets in California and provides combined strategies to maximize the revenues via unidirectional EV-Grid integrations (V1G). EV itinerary and usage information from a commercial demonstration site is utilized to model the EV flexibilities, based on which, we modeled the heterogeneous market rules using mixed-integer programming approaches. The system cost-saving performance is analyzed with respect to fleet properties and market constraints, including flexibility, participation threshold, and baseline calculaton, etc.

## I. Nomenclature

**Indices and Sets**

| | |
|---|---|
| $D, d$ | days in each planning month, day index |
| $D_{PDP}, d_{PDP}$ | days with PDP event, PDP day index |
| $T, t$ | time steps for one day, time index |
| $T_{PDP}, t_{PDP}$ | time steps during PDP event, time index in PDP days |
| $T_P$ | time steps in peak periods; |
| $T_{PP}$ | time steps in part-peak periods; |
| $T_M$ | time steps in any time maximal periods; |
| $I, i$ | demand charge periods including peak, part-peak and any time maximal periods; index of demand charge periods, i.e. 1, 2,… |
| $T_i$ | time steps for demand charge period $i$ |
| $N_p^d(t)$ | plugged-in vehicles at time $t$ on day $d$ |

**Parameters and Variables**

| | |
|---|---|
| $b_n^d(t)$ | binary charging indicator for vehicle $n$ |
| $p_n^d(t)$ | charging power for vehicle $n$ at time $t$ on day $d$ |
| $\underline{p}$ | minimal effective charging power |
| $\overline{p}$ | maximal effective charging power |
| $e_n^d(t)$ | energy charged to vehicle $n$ by time $t$ on day $d$ |
| $e_{n,req}^d$ | energy requested by vehicle $n$ on day $d$ |
| $t_n^{d,a}$ | arrival time of vehicle $n$ on day $d$ |
| $t_n^{d,f}$ | finish charging time of vehicle $n$ on day $d$ |
| $t_n^{d,l}$ | departure time of vehicle $n$ on day $d$ |
| $e_{n,d}^{+/-}(t)$ | fastest/slowest accumulated energy boundaries of vehicle $n$ at time $t$ on day $d$ |
| $E_d^{+/-}(t)$ | fastest/slowest accumulated energy boundaries of the aggregate EVs by time $t$ on day $d$ |
| $P^d(t)$ | aggregated charging power at time $t$ on day $d$ |
| $L^d(t)$ | baseload at time $t$ on day $d$ |
| $\eta_c$ | charging efficiency |
| $\lambda(t)$ | energy charge rate ($/kWh) for time $t$ |
| $\omega_i$ | demand charge rate ($/kW) for period $i$ |
| $P_{PDP}^{CR}$ | capacity reserve value for PDP policy |
| $\lambda_{PDP}$ | energy charge rate during PDP events |
| $C_{EC}$ | monthly energy charge cost |
| $C_{DC}$ | monthly demand charge cost |
| $\pi_{PDP}^p$ | PDP credit rate for peak demand period |
| $\pi_{PDP}^{pp}$ | PDP credit rate for part-peak demand period |
| $R_{PDP}^p$ | PDP credit for peak demand periods |
| $R_{PDP}^{pp}$ | PDP credit for part-peak demand periods |
| $C_{PDP}$ | monthly energy charge during PDP events |
| $R_{AS}$ | monthly revenue from ancillary service market |
| $\pi_{ru/rd}^d(t)$ | regulation up/down price at time $t$ on day $d$ |
| $P_{ru/rd}^d(t)$ | regulation-up/down bid (kW) at time $t$ on day $d$ |
| $\rho_{ru/rd}$ | utilization factor of regulation up/down signals |
| $B^d(t)$ | aggregate power baseline at time $t$ on day $d$ |
| $b_{agg}^{d,B/P}(t)$ | binary indicator for aggregate baseline/actual power at time $t$ on day $d$ |
| $b_{agg}^{d,ru/rd}(t)$ | binary indicator for aggregate regulation up-/down power at time $t$ on day $d$ |
| $b_{agg}^{d,bu/bd}(t)$ | binary indicator for aggregate power with full up/down signals at time $t$ on day $d$ |
| $\underline{P}_{ru/rd}$ | min. threshold in regulation up/down markets |
| $R_{PDR}$ | revenue in PDR markets |
| $P_{sell}^d(t)$ | virtual sell power at time $t$ on day $d$ |
| $\pi_{pdr}^d(t)$ | PDR market price at time $t$ on day $d$ |
| $b_{agg}^{sell}(t)$ | binary indicator for sell power in PDR market |
| $\underline{P}_{sell}^{PDR}$ | minimum power threshold on PDR market |
| $R_{DBP}$ | revenue from DBP market |
| $\pi_{DBP}$ | credit rate on DBP market |
| $P_{rdc}^d(t)$ | reduced power due to DBP events at time $t$ on day $d$ |
| $b_{DBP}^d(t)$ | binary indicator for DBP power reduction at time $t$ on day $d$ |

This work was funded by California Energy Commission (CEC) Contract 14-057.
Bin Wang, Rongxin Yin, and Doug Black are with Lawrence Berkeley National Laboratory, Berkeley, CA, 94720.

## II. INTRODUCTION

ELECTRIC vehicle has been identified as a valuable type of DERs to provide support for a number of grid services, including frequency regulation [1], the cost mitigation [1]–[4] and resiliency improvement [5], [6]. Vehicle-to-grid technology, usually referred to as V2G, enables the bi-directional power flow between the grid and PEVs, ranging from single vehicle to vehicle fleets [7]–[9] with larger capabilities, where the PEV batteries are aggregated [10] to provide multiple grid services.

DR is defined by California Public Utilities Commission (CPUC) [11], [12] as the strategy to reduce or change the load in response to changes in the price of electricity, financial incentives, changes in wholesale market prices, or changes in grid conditions. In retail market, Time-of-Use (TOU) tariff option is provided. For commercial customers, more DR products are offered, including demand charges on top of the energy charges in TOU tariff structure, and option to opt in peak-day-pricing (PDP) plan. In addition, the aggregated resources can also join the wholesale energy markets, i.e. proxy demand response (PDR) markets and ancillary services markets. Regulation up and down markets are representative ancillary service markets. Demand bid program (DBP) and capacity bid program (CBP) are also available for customers.

[13] explicitly models EV charging load as deferrable load, which can be shifted to different time periods in order to maximize the individual or social welfare in demand response programs. Given the global objective as a convex function, Gan [14] provided definitions of optimal charging profiles and proved that charging profiles are equivalent if there are multiple feasible solutions. He also presented decentralized and asynchronous approaches for individual vehicle control. Besides the simulation-based approaches, implementable solutions for EV control have been discussed in [2], [15]–[17], where price-based, event-based and data-driven approaches are developed and deployed with real-world systems. [7] models EVs with the integration of the ancillary service market, where day-ahead bidding is scheduled and EVs are controlled in real-time following the regulation up and down signals every 4 second. PEV travel itineraries are modeled and processed by [2], [18] in order to support grid services. Authors in [19] proposed a methodology to aggregated individual PEV batteries as one single larger virtual battery with energy and power boundaries, which is also utilized in this paper. According to the analysis in [21], battery degradation plays significant roles in V2G applications, while in V1G case, the degradation cost is negligible.

In real-world implementations, more physical and practical constraints have been observed in the demonstration project. For instance, the EV charging power cannot be continuously controlled, which confirms with the SAE J1772 [20] definitions in section 5.3.5, but contradicts with almost all existing EV scheduling models, resulting in the non-convexity of the objectives and constraints in the optimization-based approaches. This limitation has been verified by field test results.

Considering the complexity of multiple electricity markets and the real-world constraints of EVs, the authors in this paper feel it necessary to provide an updated and comprehensive modeling framework for the EV-grid integration problem. In addition, specific market rules and real-world EV charging constraints should be modeled in detail to improve the fidelity. In this paper, we select multiple California energy and DR markets to integrate EVs in a commercial demonstration site. Specifically, we consider the TOU tariff structure in PG&E territory, i.e. the E-19, including energy charge, demand charge program as the base case, which is later combined with multiple other DR products, including the peak day pricing (PDP), ancillary service market, i.e. regulation up/down markets, wholesale energy market, i.e. PDR market, and DBP markets. Integer programming approaches are intensively utilized to model the EV charger constraints and complex market rules, such as such as dis-continuity of charging power, minimal power threshold for market participation, and the minimal duration of consecutive commitment of market participations, etc. Specifically, for frequency regulation markets, we adapt the modeling and evaluation approaches in [7], inspired by the real-world regulations signals collected from California energy markets. As an evaluation approach, the main contributions of this paper can be summarized as: 1) combinative market participation strategies for EVs in California DR markets; 2) comprehensive modeling framework to support EV operations in multiple California DR markets, using integer programming approaches; 3) analysis of critical factors in the decision-making framework that have significant impact on system cost-saving performance; 4) datasets of real-world baseload and EV usage are utilized for modeling and analysis.

This paper is organized as follows: Section III provides the detailed modeling approaches in multiple California DR markets; Section IV validates the modeling framework with real-world datasets from AlCoPark Garage, a commercial site in north California with both fleet EVs and public EVs. Cost and revenue analysis are provided with the existing the EV resources, as well as the impact of the flexibility of EV fleets. Finally, section V concludes the paper with future efforts.

## III. MODELING OF EVs IN MULTIPLE ENERGY MARKET

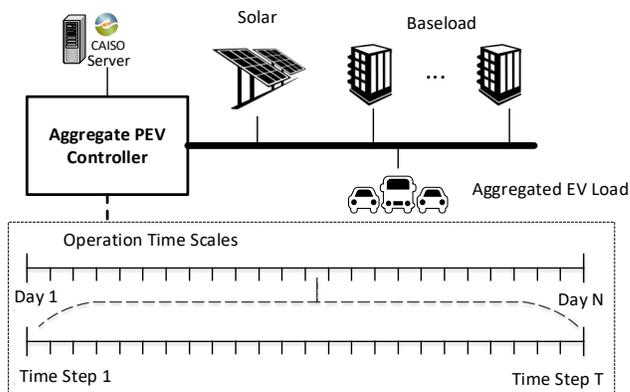

Figure 1 System overview

In this paper, EV load is modeled as deferrable load which can be shifted to different time windows to achieve various grid objectives in different energy markets. Accordingly, optimization-based strategies were developed that allow the EV fleet manager to coordinate the integration of EVs with multiple different market strategies in order to minimize the energy cost for serving the transportation required from the fleet EVs. The system architecture is shown in Figure 1. The aggregate EV controller will retrieve day-ahead pricing info from multiple DR markets from CAISO servers and collect the EV usage info, including energy demand and itineraries, from individual EV drivers. This communication has already been enabled in the demonstration project. During the next-day operation, each EV follows the day-ahead schedule in each time step to fulfill its own energy demand. For all DR markets modeled in this paper, the pricing info can be obtained in day-ahead fashion, thus no online operations are needed in this problem setting. For the real-time regulation signal, we use the utilization factors to model its impact, which is discussed in Section III.C.

*A. Aggregation of EVs*

For each individual vehicle $n$ on day $d$, the following constraints should be satisfied.

$$b_n^d(t) \cdot \underline{p} \leq p_n^d(t) \leq b_n^d(t) \cdot \bar{p} \tag{1}$$

$$e_n^d(t+1) = e_n^d(t) + p_n^d(t) \cdot \eta_c \cdot \Delta t \tag{2}$$

$$e_n^d(t_n^{d,l}) \geq e_{n,req}^d \tag{3}$$

$b_n^d(t)$ in equation (1) is the indicator of whether vehicle $n$ is charging at time $t$. Note that, the feasible charging range is not continuous so as to model the real-world EV chargers. When $b_n^d(t)$ is set to 0, both the left and right hand-sides are 0, constraining the charging power to 0, i.e. the inactive state. For the active state, the charging power threshold $\bar{p}$, i.e. minimal charging power, is set to 1.5 kW, which corresponds to the limit of the chargers used in the demonstration project. Equation (2) indicates the accrual of energy consumption for each vehicle and the energy consumption value at the time of charging session finish time $t_n^f$, i.e. $e_n^d(t_n^l)$, should be larger than the requested amount $e_{n,req}^d$. Note that energy requests for vehicles are collected by a driver-charger interface.

In order to reduce the number of decision variables in the optimization problem, modeling approaches from [19] are adapted to aggregate numerous individual EVs as one single virtual battery with power and energy boundaries, hereby improving the computational efficiency. According to this approach, any trajectory that falls between the power and energy boundaries can be achieved by controlling each PEV's charging power. The approach is summarized as follows:

$$E_d^{+/-}(t) = \sum_{n \in N_p(t)} e_{n,d}^{+/-}(t), \quad \forall t \in T \tag{4}$$

$$E_d^-(t) \leq \sum_{\tau=0}^{t} P^d(\tau) \cdot \Delta t \leq E_d^+(t), \quad \forall t \in T \tag{5}$$

$$\bar{P}^d = \sum_{n \in N_p^d(t)} \bar{p} \cdot \eta_c \tag{6}$$

$$b_{agg}^d \cdot \underline{p} \leq P(t) \leq b_{agg}^d \cdot \bar{P}^d, \quad \forall t \in T \tag{7}$$

The aggregate energy boundaries, i.e. $E_d^{+/-}(t)$, are obtained by summing up $e_{n,d}^+(t)$, which is the accumulated energy curved from the as-fast-as-possible charging pattern, and $e_{n,d}^-(t)$, which is from the as-late-as-possible charging pattern. In addition, the total power consumption value should be lower than the aggregated power from all available vehicles at time $t$. Note discontinuity of the aggregated power is also modeled, similar to that in equation (1). The optimal power consumption profiles obtain for day-ahead operations will be used as the reference for EVs to follow during the real-time operations in distributed and asynchronous fashions, which are, however, not the focus of this paper.

*B. Time-of-Use (TOU) Tariff Structure*

For commercial sites in California TOU markets, two categories of costs are generally applied to customers' bills, i.e. energy charge and demand charge. Energy charges are calculated by the product of amount of electricity, measured in kilowatt-hours (kWh) used per time period, and the per-kWh rate for those respective time periods. Demand charge is calculated by using the maximum load measurement in each demand period, multiplied by the corresponding demand charge rate, in $/kW. Thus, the total monthly cost of energy charge is modeled by equation (8), where cost of energy consumption in different time periods are all included. Equation (9) models the total monthly demand charges, where $I$ denotes the set of the demand charge periods. In the case of the E-19 tariff in the PG&E territory, there are three demand charge periods for summer months, i.e. peak, part-peak and any time max periods, while two periods in the winter, i.e. part-peak and any time max periods.

$$C_{EC} = \sum_{d \in D} \sum_{t \in T} \left( L^d(t) + \sum_{n=1}^{N_p^d(t)} P^d(t) \right) \cdot \Delta t \cdot \lambda(t) \tag{8}$$

$$C_{DC} = \sum_{T_i \in \{T_p, T_{pp}, T_M\}} \max_{t \in T_i} \left( L^d(t) + \sum_{n}^{N_p^d(t)} P^d(t) \right) \cdot \omega_i \tag{9}$$

Thus, to minimize the monthly energy bills by considering only the energy charge and demand charges, a deterministic optimization problem is formulated as:

**Problem 1 – TOU charges (energy charge + demand charge)**
Objective:    minimize $(C_{EC} + C_{DC})$
Subject to:   (1)-(9)

*C. Integration with Peak-day Pricing (PDP) plan*

In order to incentivize the electricity consumers to shift their energy consumption to periods with less system demand during specific events, peak day pricing (PDP) strategies have been applied to a number of tariff structures for both residential and commercial customers. In tariff plans, such as the E-19, the enrollment into PDP pricing on top of TOU scheme is default unless the customers select to opt out. The customer has to submit a capacity reservation value, i.e. $P_{PDP}^{CR}$, to the load serving entity, in this case, the utility company. The demand peaks in different demand periods that exceed the capacity reservation value will be protected from the demand charges by the PDP policy, i.e. credits will be billed to customers for the exceeding amount. This policy is modeled by equation (10) and (11). However, the total energy consumption in kWh below $P_{PDP}^{CR}$ during PDP events will be billed with PDP energy charge rate $\lambda_{PDP}$, which is modeled by equation (12). The optimal PEV charging problem with the PDP market participation is summarized in the problem 2.

$$R_{PDP}^{p} = \pi_{PDP}^{p} \cdot ( \max_{\substack{d \in D_{PDP} \\ t \in T_{PDP} \cap T_P}} \left( L^d(t) + \sum_{n \in N_p^d(t)} p_n^d(t) \right) - P_{PDP}^{CR}) \quad (10)$$

$$R_{PDP}^{pp} = \pi_{PDP}^{m} \cdot ( \max_{\substack{d \in D_{PDP} \\ t \in T_{PDP} \cap T_{pp}}} \left( L^d(t) + \sum_{n \in N_p^d(t)} p_n^d(t) \right) - P_{PDP}^{CR}) \quad (11)$$

$$C_{PDP} = \lambda_{PDP} \cdot \sum_{d \in D_{PDP}} \sum_{t \in T_{PDP}} \max_{\substack{d \in D_{PDP} \\ t \in T_{PDP}}} (L^d(t) + \sum_{N_p^d(t)} p_n^d(t) - P_{PDP}^{CR}, 0) \cdot \Delta t \quad (12)$$

**Problem 2 – TOU charges with PDP integration**
Objective: $minimize \ (C_{EC} + C_{DC} - R_{PDP}^{p} - R_{PDP}^{pp} + C_{PDP})$
Subject to: (1) - (12)

*D. Integration with Ancillary Service Market*

To achieve instantaneous balance between the supply and demand sides of electricity transmission system, ancillary services can be utilized by calling services from various grid components, not only the traditional electricity generators but also demand side DERs. A regulation up/down market is a representative type of ancillary service market. PEVs with the capability to follow the up and down regulation signals in a short period of time, can be coordinated to serve as effective and reliable resources to provide regulation services. Based on formulation of Problem 1, the EV integration with regulation market participation is modeled as follows:

$$R_{AS} = \sum_{d \in D} \sum_{t \in T} (P_{ru}^d(t) \cdot \pi_{ru}^d(t) + P_{rd}^d(t) \cdot \pi_{rd}^d(t)) \cdot \Delta t \quad (13)$$

$$P^d(t) = B^d(t) + \rho_{ru} \cdot P_{ru}^d(t) + \rho_{rd} \cdot P_{rd}^d(t) \quad (14)$$

$$b_{agg}^{d,B}(t) \cdot \underline{p} \leq B^d(t) \leq b_{agg}^{d,B}(t) \cdot \overline{P}^d \quad (15)$$

$$b_{agg}^{d,P}(t) \cdot \underline{p} \leq P^d(t) \leq b_{agg}^{d,P}(t) \cdot \overline{P}^d \quad (16)$$

$$E_d^-(t) \leq \sum_{t_0}^{t} P^d(\tau) \cdot \Delta t \leq E_d^+(t) \quad (17)$$

$$E_d^-(t) \leq \sum_{t_0}^{t} B^d(\tau) \cdot \Delta t \leq E_d^+(t) \quad (18)$$

$$b_{agg}^{d,rd}(t) \cdot \underline{P}_{rd} \leq P_{rd}^d(t) \leq b_{agg}^{d,rd}(t) \cdot \overline{P}^d \quad (19)$$

$$b_{agg}^{d,ru}(t) \cdot \underline{P}_{ru} \leq P_{ru}^d(t) \leq b_{agg}^{d,ru}(t) \cdot \overline{P}^d \quad (20)$$

$$b_{agg}^{d,bd}(t) \cdot \underline{p} \leq B_d(t) + P_{rd}(t) \leq b_{agg}^{d,bd}(t) \cdot \overline{P}^d \quad (21)$$

$$b_{agg}^{d,bu}(t) \cdot \underline{p} \leq B_d(t) - P_{ru}(t) \leq b_{agg}^{d,bu}(t) \cdot \overline{P}^d \quad (22)$$

Equation (13) shows the expression for calculating the total revenue from day-ahead frequency regulation markets. The revenue consists of the regulation-up capacity payment and regulation-down capacity payment. Unlike the modeling approaches in previous research where day-ahead commitments can be violated with penalties, we do not intend to violate the commitment in any circumstances due to the performance regulations in California ancillary service markets

Due to the non-continuity property of power boundaries, auxiliary binary decision variables are defined to indicate the options to participate in the regulation up and down markets. Given regulation signals from CAISO, an aggregate EV fleet will follow the signals, i.e. increase or decrease the aggregated power consumption of the EVs. However, the revenue is calculated based on the day-ahead bids, i.e. the committed regulation up and down capacities, rather than the actual increased or decreased power consumption following real-world regulation signals, indicated by equation (14). The negative (up), $\rho_{up}$, and positive (down), $\rho_{down}$, utilization factors represent the fraction of the committed regulation dispatched by the CAISO control signal. Actual utilization factors collected in a real-world demonstration project at the Los Angeles Air Force Base, were used in the simulations presented here. The baseline aggregate power $B^d(t)$, indicating the original power consumption profile assuming no regulation signals, while $P^d(t)$ is the actual power profile in equation (15) – (18). Here, $B^d(t)$ is considered as a decision variable. Equation (19), (20) model the constraints that the aggregate fleets can participate in the regulation up or down markets, or choose to stay out of the markets. We also assume that the aggregated EV fleet can follow all regulation signals, i.e. the actual power consumption should always stay in the power boundaries, which is modeled by equations (21) and (22).

Note that, the aggregator can make regulation up and down bids for the same time periods, even one of them will not be called during implementation, but still getting benefits for the bids. Additionally, the actual aggregate power and the aggregated baseline profiles should both satisfy the aggregate energy and power constraints, modeled in equation (4)- (7). The problem is formulated as:

**Problem 3 – TOU charges with regulation markets**
Objective: $minimize \ C_{EC} + C_{DC} - R_{AS}$
Subject to: (1) - (9), (13) – (22)

*E. Modeling of PDR market*

Aggregated EVs can also participate in the proxy demand resource (PDR) market, where the fleet EVs are treated as a virtual battery with flexibility to "sell" the power in the PDR market. For EVs with V2G capabilities, "sell" operations can be achieved by the discharging the vehicle batteries, while for V1G, "selling" of power would be achieved by reducing the aggregate power consumption relative to a power consumption baseline. The model is presented as follows:

$$R_{PDR} = \sum_{d \in D} \sum_{t \in T} P_{sell}^d(t) \cdot \pi_{pdr}^d(t) \cdot \Delta t \quad (23)$$

$$M_s \cdot \left(1 - b_{agg}^{PDR}(t)\right) \leq P^d(t) - B^d(t) + P_{sell}^d(t) \leq M_b \cdot (1 - b_{agg}^{PDR}(t)) \quad (24)$$

$$b_{agg}^{sell}(t) \cdot \underline{P}_{sell}^{PDR} \leq P^d(t) \leq b_{agg}^{sell}(t) \cdot \overline{P}^d \quad (25)$$

The revenue from the PDR market is a product of the virtual sell power, i.e. $P_{sell}^d(t)$, multiplied by the corresponding PDR market prices, i.e. $\pi_{pdr}^d(t)$ in equation (23). The baseline power consumption $B^d(t)$, is typically the averaged value of a number of previous days, thus here we model it as a known profile before optimization. In reality, PDR market participation has a requirement of minimal threshold for virtual sell power, i.e. $\underline{P}_{sell}^{PDR}$ in equation (25). In equation (24), $b_{agg}^{PDR}(t)$ is the binary indicator of whether the fleets are participating in the PDR market. When $b_{agg}^{PDR}(t) = 1$, i.e. participating, equation (24) is reduced to:

$$P^d(t) = B^d(t) - P_{sell}^d(t) \quad (26)$$

where the actual power consumption value $P_n^d(t)$ equals the baseline power $B^d(t)$ minus the virtual sell power $R_{sell}^d(t)$. When $b_{agg}^{PDR}(t) = 0$, indicating no participation, equation (17) evolves to

$$M_s \leq P^d(t) - B^d(t) + P_{sell}^d(t) \leq M_b \quad (27)$$

where $M_s$ is a sufficiently small number and $M_b$ is a sufficiently big number. Equation (27) remains true for all cases, making it a redundant constraint in the optimization problem, which can be effectively handled by current solvers with mixed-integer capabilities. In addition, to model the consecutive engagement constraint, some numerical approaches are applied as shown in equation (28) and equation (29),

$$b_c(t_0) = b_{agg}^{PDR}(t_0) \quad (28)$$

$$b_c(t) \leq 1 - b_{agg}^{PDR}(t - \Delta t), \quad \forall t \in T \quad (29)$$

$$b_c(t) \leq b_{agg}^{PDR}(t) \quad (30)$$

$$b_c(t) \geq b_{agg}^{PDR}(t) - b_{agg}^{PDR}(t - \Delta t), \quad \forall t \in T \quad (31)$$

$$\sum_{\tau=t}^{\min(t+N_c-1,T)} b_{agg}^{PDR}(t) - N_c \geq -M_b \cdot (1 - sc(t)), \quad (32)$$
$$\forall t \in T$$

An auxiliary binary decision variable, i.e. $b_c(t)$, is utilized to model the consecutive participation constraint. $b_c(t) = 1$ indicates the beginning of a new block of consecutive participation at $t$. Equations (28) – (32) guarantee the number of consecutive participating time steps is greater or equal to $N_c$. Incorporating binary decision variables into the optimization problems results in a mix-integer programming problem, where the sophisticated numerical solvers are needed. With PDR market integration, the overall problem is formulated as follows:

**Problem 4 - TOU charges with PDR market participation**
Objective: $\text{minimize } C_{EC} + C_{DC} - R_{PDR}$
Subject to: (1) - (9), (23) - (32)

### F. Modeling of Demand Bidding Program (DBP)

In order to increase system reliability, some utility companies are paying additional incentives to industrial, commercial or agricultural customers to reduce their energy consumption during certain time periods. An example of this, is PG&E's demand bidding program (DBP). DBP events are dispatched in day-ahead operations, so pre-planning is necessary for optimizing the benefits. The modeling approaches for DBP markets are very similar to those for PDR market, except that a fixed credit ($/kW) is used to calculate the revenues. Aggregated EVs can be utilized as valuable resources in response to DBP events following the above virtual battery modeling approaches.

$$R_{DBP} = \pi_{DBP} \cdot \sum_{d \in D} \sum_{t \in T} P_{rdc}^d(t) \cdot \Delta t \quad (33)$$

$$b_{DBP}^d(t) \cdot \underline{P}_{DBP} \leq P^d(t) \leq b_{DBP}^d(t) \cdot \overline{P} \quad (34)$$

Note that there are also consecutive participation requirements and a minimal power reduction requirement, thus the constraints for PDR market, i.e. equation (23) – (32) are also valid for DBP. The problem is defined as follows:

**Problem 5 - TOU charges with DBP market participation**
Objective: $\text{minimize } C_{EC} + C_{DC} - R_{DBP}$
Subject to: (1) - (9), (23) - (34)

## IV. RESULTS AND DISCUSSION

### A. Actual EV Charging Profiles, Electric Utility Rate, and AS Reg. Prices

Table 1 Dataset property

| Number of Sessions | Number of EVSEs | First Session Date | Last Session Date |
|---|---|---|---|
| 20,363 | 25 | 3/15/2013 | 9/7/2017 |

The real-world datasets of public and fleet EV charging at the AlCoPark Garage, including the whole building demand from the PG&E utility electric meter, from 2013 to 2017, were collected and utilized for the simulations presented and discussed below. The dataset properties are displayed in Table 1.

Table 2 PG&E E-19 demand charge and energy charge rates

| Demand Charges | $/kW | Time Period |
|---|---|---|
| Maximum Peak Demand Summer | $18.74 | 12:00 PM-6:00 PM |
| Maximum Part-Peak Demand Summer | $5.23 | 8:30 AM-12:00 PM and 06:00 PM-09:30 PM |
| Maximum Demand Summer | $17.33 | Any time |
| Maximum Part-Peak Demand Winter | $0.13 | 8:30 AM-9:30 PM |
| Maximum Demand Winter | $17.33 | Any time |
| **Energy Charges** | **$/kW** | **Time Period** |
| Peak Summer | $0.14726 | 12:00 PM-6:00 PM |
| Part-Peak Summer | $0.10714 | 8:30 AM-12:00 PM and 06:00 PM-09:30 PM |
| Off-Peak Summer | $0.08057 | Any time |
| Part-Peak Winter | $0.10166 | 8:30 AM-09:30 PM |
| Off-Peak Winter | $0.08717 | Any time |

The demonstration site is under the PG&E E-19 tariff with energy and demand rates shown in Table 2. Ancillary service regulation up and down day ahead prices were collected from CAISO's Open Access Same-time Information System (OASIS) for dates 1/1/15 to 12/31/16.

### B. Simulation results

Results of optimizations of EVs in demand response programs and ancillary service markets described above are presented here. The first example optimizes charging schedules solely to minimize electric TOU costs. The second example optimizes to minimize TOU costs and maximize ancillary service regulation revenue.

First, the load shifting and cost reduction effects of smart charging programs under only TOU prices are presented. As shown in Table *2* above, the energy charge and demand charge rates in winter are lower than those in summer. As a result, AlCoPark Garage's actual total monthly costs for energy charges in winter were slightly lower than those in summer, indicated by the blue bars in Figure *2*, and the total monthly demand charges are considerably lower than those of summer, indicated by the orange bars in Figure *2*.

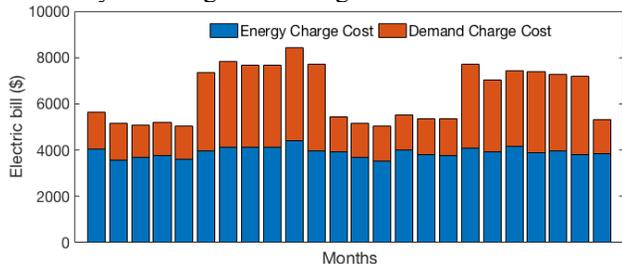

Figure 2 Total monthly electric bills from Jan. 2015 to Dec.

1) Impact of demand charge

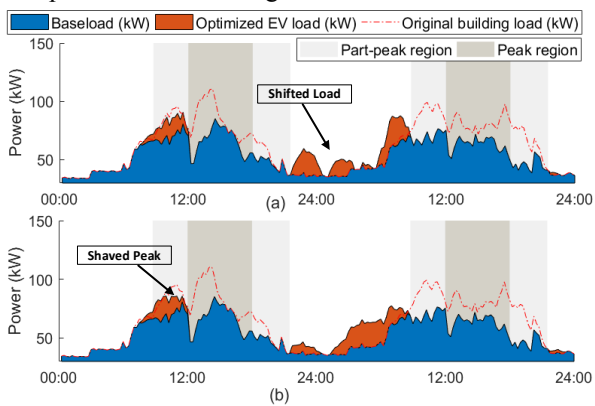

Figure 3 Load shifting of peak demand on a typical weekday

Since the demand charge is calculated by the monthly peak load of the commercial site multiplied by the corresponding demand charge rate, the smart charging program has a tendency to reduce the monthly peaks in multiple demand windows. For instance, the load profile of the day with the maximum monthly demand is shown in Figure 3. The optimized PEV load profile with only energy charge, i.e. upper Figure 3 (a), is compared with the one with both energy charge and demand charge, i.e. the lower Figure 3 (b), where the original load between 9:00AM and 12:00PM (part-peak period) is flattened.

2) Ancillary service market participation

To investigate the impact of ancillary service market integration, an additional option in the simulated smart charging program is added to allow the EV fleet to modify the aggregate power consumption profile in response to the regulation up and down prices from the CAISO ancillary services market, i.e. the problem 4 defined in the previous section. As shown in Figure 4, the EV load profile becomes spikier supporting the ancillary service market participation because the optimization tends to increase or decrease the power consumption when a high regulation up/down price is anticipated. However, with the possibility that increased EV charging load may cause to create new demand peaks, herein increasing the demand charges, the optimization has to evaluate the trade-off globally on a monthly basis. Illustrated by Figure 4, the adjusted power consumption profiles due to regulation market participation are constrained so as not to exceed the monthly demand peaks set by the TOU-based optimization.

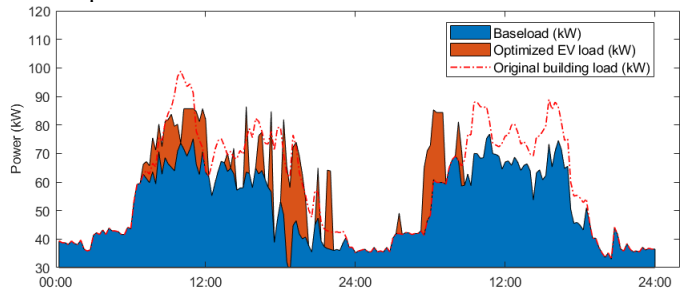

Figure 4 Example Load Profile for Two Days of Ancillary Service Regulation Up and Down Market Participation

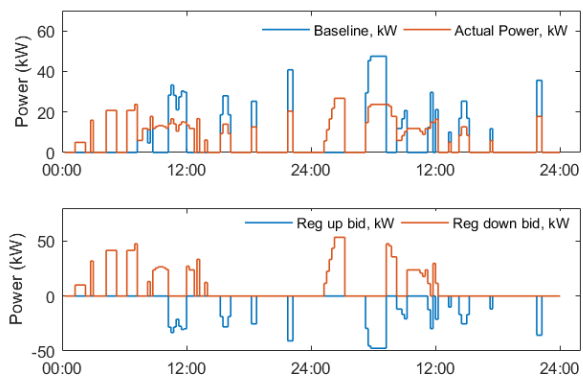

Figure 5 Details of regulation market participation

The details of the regulation market participation are shown in Figure 5, including the baseline power, actual EV load profile (upper), and the actual regulation up/down bids (lower). Note that in (a), the blue curve denotes the EV baseline load profile and the red curve is the actual EV power consumption curve. Using both curves in the optimization, energy consumption (kWh) was held constant, constrained by equations (17) and (18). Note that the duration of each regulation commitment was assumed to be 15 minutes in the optimization, within which the actual regulation signals are dispatched every 4 seconds. A 15-minute interval was also considered as the finest resolution for EV control. In addition, both regulation up and regulation down bids were allowed in the same time periods. Due to the assumption about the regulation up/down utilization rates, the regulation up/down bids were called partially, and the adjusted EV power consumption was reflected on as the differences between the baseline (blue) and the actual load profile (red) in Figure 5 (upper).

The monthly revenue results were collected by simulating EV management strategies for each month from Jan. 2015 to Dec. 2016 and shown in Table 3. The highest monthly revenue was $115 in Dec. 2016 and the lowest revenue was $78 in Sep. 2016. The relationship between monthly regulation revenue and EV charging flexibility at the AlCoPark Garage is shown in Figure 6. To represent the monthly average distance between the upper and lower boundaries of the virtual battery, the flexibility index of the aggregated EVs is defined as:

$$f_{agg} = \frac{\sum_{d \in D}\sum_{t \in T}(E^+(t) - E^-(t))}{D \cdot T} \qquad (35)$$

As indicated in Figure 5, the ability to generate profits in regulation markets is positively correlated with the flexibility index of the aggregated virtual battery, with a correlation coefficient of 0.667.

Table 3 Monthly revenue from regulation market

| Year | Month | Revenue | Year | Month | Revenue |
|---|---|---|---|---|---|
| 2015 | 1 | $90.38 | 2016 | 1 | $74.53 |
| 2015 | 2 | $66.42 | 2016 | 2 | $66.33 |
| 2015 | 3 | $73.88 | 2016 | 3 | $90.37 |
| 2015 | 4 | $95.73 | 2016 | 4 | $71.75 |
| 2015 | 5 | $78.09 | 2016 | 5 | $68.62 |
| 2015 | 6 | $67.75 | 2016 | 6 | $95.70 |
| 2015 | 7 | $76.57 | 2016 | 7 | $69.22 |
| 2015 | 8 | $68.51 | 2016 | 8 | $77.28 |
| 2015 | 9 | $78.16 | 2016 | 9 | $61.50 |
| 2015 | 10 | $80.48 | 2016 | 10 | $78.24 |
| 2015 | 11 | $63.24 | 2016 | 11 | $79.70 |
| 2015 | 12 | $98.11 | 2016 | 12 | $115.14 |

3) PDR market participation

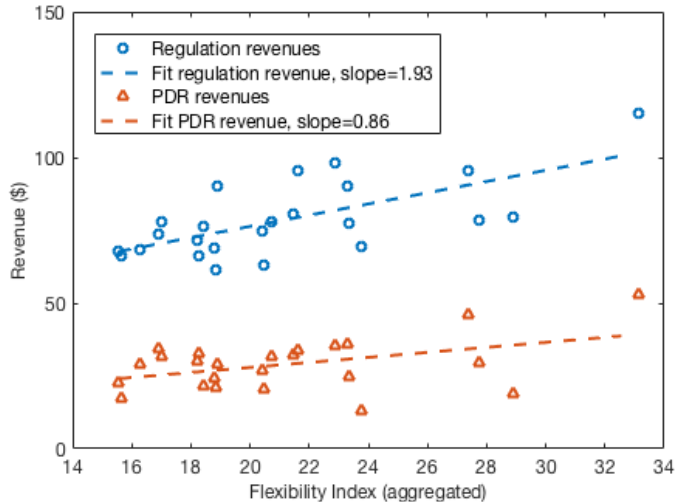

Figure 6 Monthly profits vs. flexibility index

Problem 5 was addressed to simulate PDR market participation. CAISO requires that each commitment into the PDR market have at least a one-hour commitment. PDR market commitments of 1-hour were modeled with the constraints represented in equations (28) – (32). As shown in Figure 7 (lower), the green curve indicates the actual EV power consumption profile, while the red curve represents the virtual sell power of the aggregated EVs given price signals from the PDR market. Note that the total energy consumption value following the actual power consumption profile should be equal to the one that follows the baseline profile generated by problem 2. In addition, problem 5 models the opportunities of the EVs to participate in the PDR market as discrete options, i.e. the EV aggregator does not have to stay in the market for the whole day and it can plan to step out of market when the PDR prices are not optimal.

The actual monthly revenues from PDR markets illustrated by the red triangles in Figure 6, where the varying flexibilities of EV fleets to generate profits from PDR markets are shown. Note the consecutive commitment constraint is set to 1 hour for the PDR market optimizations.

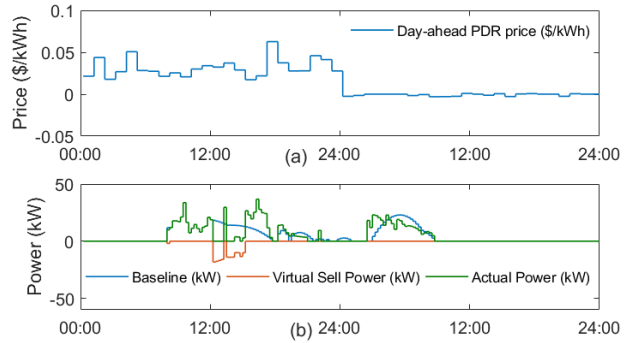

Figure 7 PDR market participation

4) DBP Participation

Table 4 Monthly revenue from DBP market

| Year | Month | Event # | Revenue |
|---|---|---|---|
| 2016 | 6 | 5 | $16 |
| 2016 | 7 | 6 | $10 |
| 2016 | 8 | 2 | $0 |
| 2016 | 9 | 1 | $0 |

The modeling approaches for the PDR market integration were similar to those for the DBP market, however, there were different requirements for commitments in the DBP market. For instance, participation in the DBP market only occurs when the DBP events are issued by the program facilitator, PG&E, however, hourly price signals in the PDR market were available on a daily basis. In PG&E's DBP program, $0.5 per kW is credited to commercial customers when they reduce their demand during DBP events. It was assumed, in this analysis, that the threshold to participate was greater than or equal to 10 kW and each commitment had to have a duration at least 2 consecutive hours. The simulation results for the DBP market participation is shown in Table 4. Due to the 2-hour commitment constraint, the existing EV resources were not qualified to participate in all of the DBP events in 2016. Thus, the profit-generating capacity for EVs is not as high as that for the regulation market, considering the limits of the fleet size and V1G power and flexibility.

5) PDP Participation

Participating the PDP program, the annual electric bill savings were expected to improve as the monthly peak demand and part-peak demand were partially protected by the Capacity Reserve Level (CRL), which was required for PDP program enrollment. Specifically, as modeled by equations (10) – (12), the monthly peaks above the CRL received PDP credits, while the energy usage not protected by CRL was billed with a fixed PDP rate. PDP events were only issued during summers, and only peak and part-peak demands were considered. The monthly PDP benefit was calculated as $R_{PDP}^{p} + R_{PDP}^{pp} - C_{PDP}$. As shown in Figure 8, PDP benefits for summer months in the year of 2016 were computed with varying CRLs. For months with only 1 PDP event, i.e. August and September, PDP credits dominated the total benefit, which decreased as the CRL increased. On the contrary, event energy charge became dominant in months with more PDP events since there was less unprotected energy usage as the CRL increased between 10 kW and 60 kW. As CRL increased greater than 60 kW, the month benefits decreased because of the weaker protection by CRL. In addition, the annual total PDP benefit varied with the CRL with the optimal CRL value close to 40 kW.

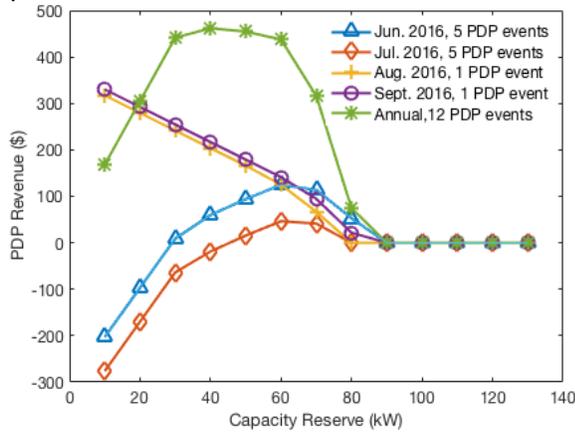

Figure 8 Impact of capacity reserve on PDP benefits

*C. Discussion*

Impacts of different factors on the revenue-generating capability of the EV fleet, including the freedom of baseline power profile selection, flexibility of individual EV, and market participation threshold, etc. will be discussed here. Simulation results indicate that proper tuning of these factors can lead to significant improvements in revenue generation.

1) Impact of baseline calculation

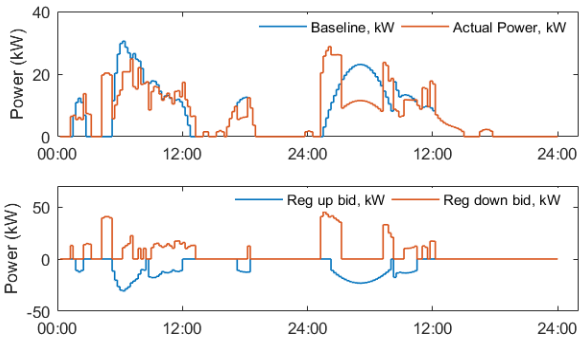

Figure 9 Fixed baseline case

As opposed to the regulation simulation described above where baseline charging power was a decision variable, here, instead, the charging profiles obtained by solving problem 2 were used as baselines. As shown in Figure 9, the actual power (red) generally follows the baseline power profile (blue), unlike that shown in Figure 5. However, the capability of the smart charging program was limited in exploring more space to generate revenues from regulation markets. The monthly revenue from regulation in June, 2016 was reduced from $95 to $36. Thus, with a preset baseline power profile, the flexibility was limited as well as the revenue-generating capability.

2) Impact of the flexibility on regulation revenue

The impact of individual vehicle charging flexibility on regulation market revenues is examined here. The total connected duration of each EV was increased by multiple ratios to simulate different degrees of EV connected time flexibilities. For the months shown in Figure 10, January, June and August in 2016, the revenues increased rapidly as the ratio of connected time to charging time increased from 1 to 2.5. However, as the ratio increased beyond 2.5, the total monthly revenue gradually plateaus. Note that in the simulation, 2-days was the maximum connected duration for each EV. Thus, within the given time period, there exist limitations of revenue improvement by extending EV connected duration flexibility. However, for real-world operations, it will be beneficial for an EV fleet manager to maximize each EV's connected duration flexibility, by either having them arrive earlier or leave later.

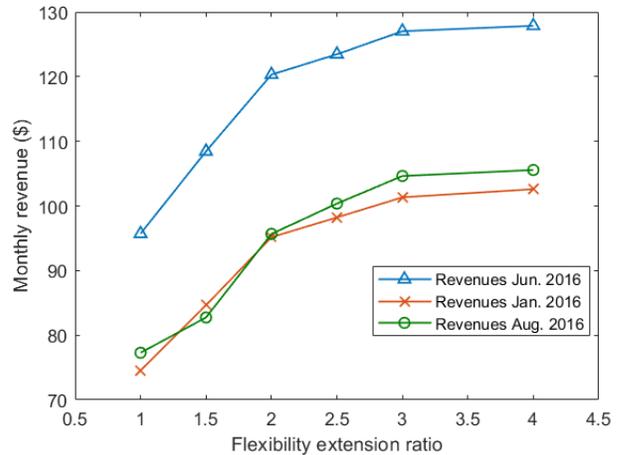

Figure 10 Impact of flexibility

3) Impact of minimum participation threshold

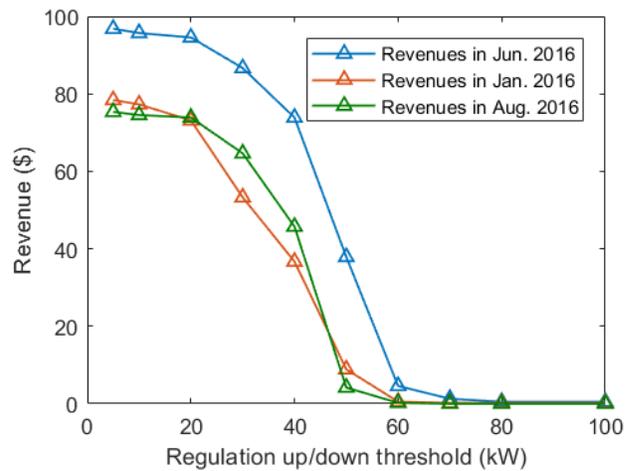

Figure 11 Impact of regulation threshold

In the market participation simulations presented above, the minimum threshold to participate was 10 kW, which was appropriate for the size of the EV fleet in this study. Simulations to investigate the impact of varying threshold values on revenue is presented here. As the threshold value increases in Figure 11, the initial revenue drop is small, however, revenue sharply decreases as the threshold increases from 20 kW to 60 kW, indicating most of the commitments failed to satisfy the constraints defined by equation (19) and (20) because the required power adjustments exceeded the capacity of the EV fleet.

## V. CONCLUSION AND FUTURE WORK

We provide comprehensive evaluation models on the integration of PEVs into various California demand response markets, analyze the revenues from them and investigate the impact of multiple factors on the revenue-generating capabilities. For the future work, we will explore the online strategies within this framework as potential implementable solutions and extend the developed strategies for V2G cases..

## VI. ACKNOWLEDGMENT

This work was funded by California Energy Commission (CEC) Contract 14-057.

## VIII. BIOGRAPHIES

**Bin Wang** is currently a postdoctoral research fellow at Lawrence Berkeley National Laboratory. His research interests include optimization and control strategies for intelligent transportation system and sustainable energy system. Bin received his Bachelor's degree in vehicle engineering from Jilin University, China in 2012, and his Ph.D. degree in mechanical engineering from University of California, Los Angeles in 2016.

**Rongxin Yin** is a Senior Scientific Engineering Associate in the Grid Integration Group of Lawrence Berkeley National Laboratory. He has been at LBNL since 2007 and focuses his career on building science, energy and buildings, and building to grid. His research focuses on energy efficiency and demand response (EE&DR) in commercial buildings, energy modeling, methodologies for baseline models for EE & DR, thermal energy storage, and smart grids. He is the primary developer of the Demand Response Quick Assessment Tool (DRQAT) for commercial buildings and refrigerated warehouses. He has a Master's of Science Degree in Building Science from the University of California, Berkeley, and a Master's of Science Degree in Mechanical Engineering from Tongji University in China.

**Doug Black** is a Mechanical Engineer in the Energy Storage and Distributed Resources Department. He is a major contributor to LBNL's healthcare energy efficiency efforts including the development of hospital benchmarking metrics and procedures and measurements of medical equipment and miscellaneous electric loads in hospitals. He is also co-leading a study to create an automated method for tuning models of building energy consumption (e.g. EnergyPlus models) using measured sub-meter data. His EnergyPlus experience includes the creation and calibration of highly resolved models of commercial buildings as well as parametric studies using prototype building models. He has a BS in Electrical Engineering from the University of Michigan and an MS and Ph.D. in Civil and Environmental Engineering from the University of California at Berkeley.